\documentclass{article}
\usepackage{amsthm}
\usepackage{amsmath,amssymb}

\newtheorem{theorem}{Theorem}

\newtheorem{definition}[theorem]{Definition}
\newtheorem{parrafo}[theorem]{}

\numberwithin{theorem}{section}
\numberwithin{equation}{theorem}

\DeclareMathOperator{\Max}{\underline{Max}}

\DeclareMathOperator{\Fnu}{\boldsymbol{\nu}}
\DeclareMathOperator{\Reg}{Reg}
\DeclareMathOperator{\Sing}{Sing}

\def\O{{\mathcal{O}}}

\title{A  Proof of Desingularization over fields of characteristic zero.}

\author{Encinas, Santiago  \\
Dep. Matem\'atica Aplicada Fundamental \\
Universidad de Valladolid \\
\texttt{sencinas@maf.uva.es}
\and
Villamayor, Orlando  \\
Dep. Matem\'aticas \\
Universidad Aut\'onoma de Madrid \\
\texttt{villamayor@uam.es}
}
\date{September 2001}

\begin{document}
\maketitle
\begin{abstract}
	We present a proof of embedded desingularization for closed ÊÊÊÊÊ
	subschemes which does not make use of Hilbert-Samuel function and
	avoids Hironaka's notion of normal flatness.
	
	This proof, already sketched in [A course on constructive
	desingularization and equivariance.  In {\emÊÊÊÊÊÊÊ ÊÊÊÊÊ
	Resolution of singularities (Obergurgl, 1997)}, vol.  181 {\em
	Progr.  Math.}, Birkh\"auser, 2000.]  page 224, is done by
	showing that desingularization of a closed subscheme $X$, in
	a smooth sheme W, is achieved by taking an algorithmic
	principalization for the ideal $I(X)$, associated to the
	embedded scheme $X$.
\end{abstract}



\section*{Introduction.}

In his monumental work (\cite{Hironaka1964}), Hironaka proved
desingularization, and also a strong form of principalization, or
say monomialization,  of ideals in regular schemes; both results
proved over fields of characteristic zero. His proof is
existential in the sense that it does not provide an algorithm in
order to achieve desingularization.

Different constructive proofs, following Hironaka's approach,
appear in \cite{BierstoneMilman1997},
\cite{EncinasVillamayor1998}, and in \cite{Villamayor1989}(see
also \cite{Villamayor1992}). Each one of these proofs provides an
algorithm of desingularization which indicates where to blow-up in
order to eliminate the singularities in a step by step procedure.
The idea is to define invariants of singular points, and to show
that these invariants improve when blowing up the set of worst
singularities. Desingularization is then achieved by repeatedly
blowing up the set of worst points.

All algorithmic procedures mentioned above make use, as Hironaka
did in his original work, of the Hilbert Samuel function. Namely
the invariant attached to a singular point, consists of the full
Hilbert Samuel function at the given point, followed by other
data.

The purpose of this work is to show that embedded desingularization
can be achieved in a different way, where the Hilbert Samuel functions
are avoided.  In fact we show that desingularization follows from
algorithmic principalization (see also \cite{EncinasVillamayor2000}
page 224).

Here principalization of ideals is meant in a strong form, in
which the ideal becomes locally monomial after a suitable sequence
of monoidal transformations. Algorithmic principalization of
ideals, is stated here in \ref{AlgResol}  as an algorithm of
resolution of \emph{basic objects} (see also
\ref{monomialization}).

This paper is organized so that the reader can easily get into
what is new in this proof, without having to go through
technicalities. To this end we do not include here details of
algorithmic resolution of basic objects; here we  focus on showing
why this proof arises quite directly from properties extracted
from the algorithm. Let us mention that the algorithm of
resolution of basic objects treated in
\cite{EncinasVillamayor2000}, and implemented in MAPLE by Bodn\'ar
and Schicho, is available at
\begin{center}
http: //www.risc.uni-linz.ac.at/projects/basic/adjoints/blowup
\end{center}
and we encourage the reader to test on examples.

This proof has also led to stronger formulations of
desingularization (see \cite{BV1} and \cite{EH}), and has been
applied in the study of desingularization of families of embedded
schemes (see \cite{ENV}).

In this paper we address equivariant desingularization for
equidimensional schemes which are  globally embedded in smooth
schemes. We finally refer to \cite{BEV} for an extension of this
proof to the class of locally embedded excellent schemes; and for
a full description of this equivariant algorithm in this more
ample context. That paper also contains a detailed discussion of
the result in \cite{BV2}, and applications.

\section{Embedded desingularization.}

\begin{definition}
\label{normalc} Let $W$ be a regular scheme and let
$Y_1,\ldots,Y_k\subset W$ be a set of closed subschemes.  We say
that $Y_1\cup \ldots\cup Y_k$ have {\em normal crossing at a point
$\xi\in W$}  if there exists a regular system of parameters
$\{x_1,\ldots,x_d\}\subset {\mathcal O}_{W,\xi}$, such that for
each $i\in \{1,\ldots,k\}$, either $\mathcal{I}(Y_i)_{\xi}={\mathcal
O}_{W,\xi}$,  or $$\mathcal{I}(Y_i)_{\xi}=\langle
x_{i_1},\ldots,x_{i_{s_i}}\rangle$$ for some  $
x_{i_1},\ldots,x_{i_{s_i}}\in \{x_1,\ldots,x_d\}$.  We say that
$Y_1\cup\ldots\cup Y_k$ have {\em normal crossing in $W$}  if they
have normal crossing at any point of $W$.
\end{definition}

\begin{definition}
\label{DefPair} Let \( W \)  be  a pure dimensional scheme, smooth
over a field \( \mathbf{k} \) of characteristic zero, and let \(
E=\{H_{1},\ldots,H_{r}\} \) be a set of smooth hypersurfaces in \(
W \) with  normal crossing. The couple  \( (W,E) \) is said to be
a {\em pair}.
\end{definition}

\begin{parrafo} \textbf{Transformation of pairs.}
A regular closed subscheme \( Y\subset W \) is said to be
\emph{permissible} for a pair \( (W,E) \) if \( Y \) has normal
crossing with  \( E \) (i.~e. with \( \cup_{i=1}^{r} H_{i} \)).

Given \( (W,E) \) and \( Y \) as above, let
\[ W \longleftarrow W_{1} \]
be the blow up at \( Y \), and set \(
E_{1}=\{H'_{1},\ldots,H'_{r},H'_{r+1}\} \), where \( H'_{i} \)
denotes the strict transform of \( H_{i} \), and \(
H'_{r+1}=\Pi^{-1}(Y) \) is the exceptional hypersurface in \(
W_{1} \). The permissibility of $Y$ insures that \(W_{1}\) is
smooth, and that \(E_{1}\) has normal crossing. We say that
\[ (W,E)\longleftarrow(W_{1},E_{1}) \]
is a \emph{transformation} of pairs defined by the permissible
center $Y$.

A sequences of transformations of pairs:
\begin{equation}
     (W_{0},E_{0})\longleftarrow(W_{1},E_{1})\longleftarrow
         \cdots\longleftarrow (W_{k},E_{k})
     \label{SeqPairs}
\end{equation}
with centers \( Y_{i} \), \( i=0,1,\ldots,k-1 \) is a composition
of transformations.
\end{parrafo}

\begin{definition} \label{equivPair}
      We say that an isomorphism \( \Theta:W_{0}\to W_{0}
      \) defines an isomorphism on the pair \( (W_{0},E_{0}) \) if \(
\Theta(H_i)\subset H_i
      \) for any \( H_i \in E_{0} \).

      A group \(G\) is said to act on a pair \((W_{0},E_{0})\), if it
acts on $W_0$, and if any \(
      \Theta\in G \) defines an isomorphism on the pair.

Let now
\[ (W_0,E_0)\longleftarrow(W_{1},E_{1}) \]
be a transformation of pairs defined by a permissible center $Y
\subset W_0$. Assume that a group $G$ acts both on the pair
$(W_0,E_0)$ and also on $Y \subset W_0$. In this case $G$ also
acts on the pair $(W_{1},E_{1})$ (\cite[Lemma 4.2]{BEV}).

      A sequence of transformations of pairs (\ref{SeqPairs})
      is said to be \( G \)-equivariant if $ G$ acts on the pairs
$(W_{i},E_{i}) $, $0 \leq i \leq k-1 $, and
     \( \Theta(Y_{i})=Y_{i} \) for $0 \leq i \leq k-1 $,
     for any  \( \Theta \in G \).
In such case the group $G$ acts on the pair $(W_{k},E_{k})$.

Note that by a step by step lifting of the action, we can also
state that (\ref{SeqPairs}) is equivariant if $G$ acts on acts on
the pair $(W_{0},E_{0}) $, and also on the centers \( Y_{i} \) for
$0 \leq i \leq k-1 $.

\end{definition}

\begin{theorem} \label{MainTh}
      \textbf{Embedded Desingularization.}

      Given a closed reduced and equidimensional subscheme
      \(X_0\subset W_0\),
      there is a sequence of transformations of pairs
      \[ (W_0,E_0=\emptyset)\longleftarrow \cdots \longleftarrow (W_r,E_r) \]
      inducing a proper birational morphism
      \( \Pi:W_{r}\longrightarrow W_{0} \), so that
      setting $Sing(X_0)$ as the singular locus of $X_0$, $Reg(X_0)=
X_0- Sing(X_0)$, and
       \( X_{r}\subset W_{r} \) the strict transform of \(X_0\), then:
      \begin{description}
     \item[(i)]  The morphism \( \Pi \) defines an isomorphism
     \[ W_{0}\setminus\Sing(X_{0})\cong W_{r}\setminus
     \bigcup_{H\in E_{r}}H \]
     and hence
     \( \Reg(X)\cong\Pi^{-1}(\Reg(X))\subset X_{r} \) via \( \Pi \).

     \item[(ii)]  \( X_{r} \) is regular and has normal crossing
     with \( E_{r} \).

     \item[(iii)] (Equivariance) Any action
     of a group \(G\) on \( X_{0}\subset W_{0} \) has a unique
     natural lifting to an action on \((W_r,E_r)\) and on
     \(X_r\subset W_r\)(\ref{gaobo}).
      \end{description}
\end{theorem}

\begin{proof}
     See \ref{ProofMain}.
\end{proof}

\section{Basic Objects.}

\begin{parrafo}
  We will prove Theorem \ref{MainTh} as a corollary of the
      algorithm of monomialization. We will essentially unify
      monomialization and desingularization by means of the notions
      of basic objects and of algorithmic resolution of basic objects,
which we now introduce.

\end{parrafo}

\begin{definition} \label{DefBasic}
      \cite[definition 1.2]{EncinasVillamayor1998}
      A \emph{basic object} \( (W_{0},(J_{0},b),E_{0}) \), is a pair \(
      (W_{0},E_{0}) \), an ideal \( J_{0}\subset\O_{W_{0}} \), and a
      positive integer \( b \). We require that \( (J_{0})_{\xi}\neq 0
      \) for any \( \xi\in W_{0} \).

  Note here that $\O_{W_{0},\xi}$ is a local regular ring; let
$\Fnu_{J_{0}}(\xi)$ denote the order of $(J_{0})_{\xi}$ at
the local ring $\O_{W_{0},\xi}$ (the biggest integer such that the
corresponding power of the maximal ideal contains
$(J_{0})_{\xi}$).

We finally define
      \[ \Sing(J_{0},b)=\{\xi\in W_{0}\mid \Fnu_{J_{0}}(\xi)\geq b\}
     ( \subset W_{0}) \]
      which is a closed subset in $W_0$.
\end{definition}

\begin{parrafo}
      \textbf{Transformation of basic objects.} \cite[definition
1.4]{EncinasVillamayor1998}\label{tbo}
      We shall say that \(
      Y_{0} \) is \emph{permissible} for the basic object \(
      (W_{0},(J_{0},b),E_{0}) \) if \( Y_{0} \)
      is permissible for the pair \( (W_{0},E_{0}) \), and, in addition,
      \( Y_{0}\subset\Sing(J_{0},b) \).  In such case let \( W_{0}\longleftarrow
      W_{1} \) be the blow-up with center \( Y_{0} \) and denote by \(
      H_{1}\subset W_{1} \) the exceptional hypersurface.  In the
      particular case in which \( Y_{0} \) is irreducible, with
generic point $\xi$, then $\xi \in \Sing(J_{0},b)$ and
$\Fnu_{J_{0}}(\xi)\geq b$ (since  \( Y_{0}\subset\Sing(J_{0},b))
\).

Set $c_{1}=\Fnu_{J_{0}}(\xi)$; there is a factorization
     $$ J_{0}\O_{W_{1}}=I(H_{1})^{c_{1}}\bar{J}_{1} $$
for a well defined sheaf of ideals  $\bar{J}_{1} \subset
\O_{W_{1}}$. In the general case, in which \( Y_{0}\) is not
necessarily irreducible,
      we obtain, in a similar way, a well defined expression as above, where now
      \( c_{1}\geq b\) is locally constant on \( H_{1} \).

We define
      $$J_{1}=I(H_{1})^{c_{1}-b}\bar{J}_{1} $$
  and set
      \[ (W_{0},(J_{0},b),E_{0})\longleftarrow (W_{1},(J_{1},b),E_{1}) \]
      which we call a \emph{transformation} of basic objects.

It should be noted, that in general, the sheaf of ideals $J_1$ is
not the {\em strict transform} of $J_0$.
\end{parrafo}

\begin{definition} \label{defResol}
      A sequence of transformations of basic objects
      \begin{equation}
         (W_{0},(J_{0},b),E_{0}) \longleftarrow\cdots\longleftarrow
         (W_{k},(J_{k},b),E_{k})
         \label{SeqBasObj}
      \end{equation}
is a \emph{resolution} of
      \( (W_{0},(J_{0},b),E_{0}) \) if \( \Sing(J_{k},b)=\emptyset \).
\end{definition}

\begin{parrafo}\label{monomialization}
Assume, for simplicity, that $E_0=\emptyset$ so that $ \cup H_i$
($H_i \in E_{k}$ in \ref{SeqBasObj}) is the exceptional locus of
the composition $ W_0 \longleftarrow W_k $. It follows from the
notion of transformation of basic objects, that if
(\ref{SeqBasObj}) is a resolution then:

  $$ J_0 \O_{W_k}= J_k . \mathcal{M} $$ where $\mathcal{M} $ is an
invertible sheaf of ideals supported on $ \cup H_i$ (locally
spanned by a monomial), and $J_k$ is a sheaf of ideals with order
at most $b-1$ at points of $W_k$.

In particular, if $b=1$ in the resolution \ref{SeqBasObj},
$J_k=\O_{W_k}$ so $$ J_0 \O_{W_k}=  \mathcal{M} $$ Hence, the
total transform of $J_0$ is an invertible sheaf of ideals, locally
defined by a monomial (monomialization).
\end{parrafo}

\section{Equivariance.}

\begin{parrafo}
      \label{gaobo}
Let $G$ be a group acting on $W$, and $X \subset W$ a subscheme;
we say that $G$ acts on $X \subset W$ when the action on $W$
induces, by restriction,  an action on $X$. If $X \subset W$ is
simply a closed set we view it as a subscheme with the unique
reduced structure.

Since we will prove both desingularization and monomialization by
means of resolutions of basic objects, we discuss now the notion
of group actions within that context.
\end{parrafo}
\begin{definition} \textbf{Group actions on basic objects.} \label{defequi}
      Consider a basic object $(W_{0},(J_{0},b),E_{0})$ and a group
$G$ acting on the pair $(W_{0},E_{0})$.
      We will say that $G$ acts on the basic object when following
conditions holds:

0) The group $G$ acts on $\Sing(J_{0},b) \subset W_0$, namely $$
\Theta(\Sing(J_{0},b))=\Sing(J_{0},b) $$ for any $\Theta \in G$.

k) Whenever a sequence of transformations of basic objects
\begin{equation}
         (W_{0},(J_{0},b),E_{0}) \longleftarrow\cdots\longleftarrow
         (W_{k},(J_{k},b),E_{k})
         \label{SeqBasObj1}
      \end{equation}
is such that the induced sequence of pairs

\begin{equation}
         (W_{0},E_{0}) \longleftarrow\cdots\longleftarrow
         (W_{k},E_{k})
      \end{equation}
is $G$-equivariant (\ref{equivPair}), then the group $G$ acts on
$\Sing(J_{k},b) \subset W_k$, namely $$
\Theta(\Sing(J_{k},b))=\Sing(J_{k},b) $$ for any $\Theta \in G$.
\end{definition}

\begin{parrafo} \textbf{Main Example.}\label{MainEx}
  Note that the previous definition involves \emph{all} possible
$G$-equivariant sequences of transformations of basic objects.

If a group $G$ acts on $W_0$, each $ \Theta \in G$ defines an
isomorphism $\Theta^{\#}: \O_{W_{0}} \to \O_{W_{0}}$. Suppose now
that a group $G$ acts on a pair \( (W_{0},E_{0}) \), and that
$J_0$ is a is $G$-invariant sheaf of ideals (i.e. \(
\Theta^{\#}(J_{0})=J_{0} \) for any $\Theta \in G$). We claim that
in these conditions, for any $b$, the group $G$ acts on the basic
object \( (W_{0},(J_{0},b),E_{0}) \).

  So assume that \(
      Y_{0} \) is permissible for the basic object \(
      (W_{0},(J_{0},b),E_{0}) \), and $G$-invariant. Note that \(
Y_{0}\subset\Sing(J_{0},b) \)
      is permissible for the pair \( (W_{0},E_{0}) \). Set \(
W_{0}\longleftarrow
      W_{1} \) the blow-up with center \( Y_{0} \), and denote by \(
      H_{1}\subset W_{1} \) the exceptional hypersurface.

Since $Y_0$ is closed and regular, it is the disjoint union of
irreducible components $Y_0= Z_1 \cup...\cup Z_s$, each $Z_i$ with
generic point, say $\xi_i \in \Sing(W_0,b)$. Note that $H_1$ is
also a union of $s$ irreducible components, say $H_1=V_1
\cup...\cup V_s$.

Let $$ J_{0}\O_{W_{1}}=I(H_{1})^{c_{1}}\bar{J}_{1} $$ and
  $$J_{1}=I(H_{1})^{c_{1}-b}\bar{J}_{1} $$
  be as in (\ref{tbo}), and note that the locally constant function
$c_{1}$ is constant and
  equal to $\Fnu_{J_{0}}(\xi_i)$ along each component $V_i$.

On the other hand, the group $G$ acts on the set $\{
\xi_1,...,\xi_s \}$,
  and if $\Theta(\xi_i)=\xi_j$ for some $\Theta \in G$, then
$\Fnu_{J_{0}}(\xi_i)=\Fnu_{J_{0}}(\xi_j)$, since $\Theta
(J_0)=J_0$.

Since the action of $G$ can be lifted to $W_1$, and since $J_0$ is
$G$-invariant, it follows that the total transform $
J_{0}\O_{W_{1}}$ is $G$-invariant.

We leave it to the reader to check that the sheaves of ideals
$I(H_{1})^{c_{1}}$, $I(H_{1})^{c_{1}-b}$, and $J_{1}$ are
$G$-invariant in $W_1$. In particular $G$ acts on $\Sing(J_{1},b)
\subset W_1 $. A step by step argument shows now that $G$ acts on
the basic object $(W_{0},(J_{0},b),E_{0})$ in the sense of
\ref{defequi}.

  Take for instance a group $G$ acting on $X\subset
W$ (\ref{gaobo}), where $X$ is a subscheme defined by a sheaf of
ideals $J$. In this case $J$ is $G$-invariant, so $G$ acts on the
basic object $(W, (J,1),\emptyset)$.

\end{parrafo}

\section{Algorithms of Desingularization.}

\begin{parrafo} We now discuss the notion and properties of algorithmic
resolution of basic object, which will lead us to a constructive
proof of the theorem of desingularization.

\end{parrafo}

\begin{definition}
     Fix a  totally ordered set \( (I,\leq) \), and a closed set \(
F\subset W \). A function \[ h:F\longrightarrow I \] is said
      to be upper-semi-continuous if:

i)  \( h \) takes only finitely many values,
      and

ii) \( \{\xi\in F\mid h(\xi)\geq \alpha\} \) is a closed set, for
      any \( \alpha\in I \).

      We denote by $$ \max{h} $$ the maximum value in $I$ achieved by
\( h \), and set $$ \Max{h}=\{\xi\in F\mid h(\xi)=\max{h}\} $$
which is closed in $F$.

      Given a basic object \( (W,(J,b),E) \), we say that an
upper-semi-continuous
      function \[ h:\Sing(J,b)\longrightarrow I \] is
      \emph{equivariant} if, for any group $G$ acting on this basic object,
      \[ h(\xi)=h(\Theta(\xi)) \qquad \forall\ \xi\in\Sing(J,b)
\qquad \forall\ \Theta\in G.\]
\end{definition}

\begin{parrafo}\label{comequiv}
Note that if \( h \) is equivariant, then \emph{any} group \( G \)
acting on the basic object \( (W,(J,b),E) \) also acts on the
closed set \( \Max{h} \subset W\). If, in addition, $\Max{h}$ is
permissible for $(W,E)$, then any such $G$ also acts on the
transform of the basic object with center $\Max{h}$.
\end{parrafo}

\begin{parrafo} \label{AlgResol}
     \textbf{Algorithm of resolution of basic objects.}
     Let \( d \) be a non-negative integer.
     An algorithm of resolution for \( d \)-dimensional basic objects
     consists of:
\end{parrafo}
\begin{description}
     \item[A]  A totally ordered set \( (I_{d},\leq) \).

     \item[B]  For each basic object \( (W_{0},(J_{0},b),E_{0}) \) with
     \( d=\dim{W_{0}} \):
     \begin{description}
    \item[i] An equivariant function \(
    f_{0}^{d}:\Sing(J_{0},b)\longrightarrow I_{d} \) is defined,
    and this function has the property that \( \Max{f_{0}^{d}}
    \subset \Sing(J_0,b) \) is permissible for \(
    (W_{0},(J_{0},b),E_{0}) \).  \smallskip

         Suppose, by induction, that an equivariant sequence with centers
         \( Y_{i}\subset \Sing(J_{i},b) \), \( i=0,\ldots,r-1 \):
         \begin{equation} \label{eqMainSeqResol}
        (W_{0},(J_{0},b),E_{0})\longleftarrow\cdots \longleftarrow
        (W_{r-1},(J_{r-1},b),E_{r-1})\longleftarrow
        (W_{r},(J_{r},b),E_{r})
    \end{equation}
         together with equivariant functions
    \( f_{i}^{d}:\Sing(J_{i},b)\longrightarrow I_{d} \),
    \( i=0,\ldots,r-1 \) have been defined, and that
         \( Y_{i}=\Max{f_{i}^{d}} \). Then:

    \item[ii] If \( \Sing(J_{r},b)\neq\emptyset \), an equivariant
    function \( f_{r}^{d}:\Sing(J_{r},b)\longrightarrow I_{d} \)
    is defined, and \( \Max{f_{r}^{d}} \) is permissible for \(
    (W_{r},(J_{r},b),E_{r}) \).

    Note that \textbf{B(ii)} says that whenever \(
    \Sing(J_{r},b)\neq\emptyset \) there is an equivariant
    enlargement of \ref{eqMainSeqResol} with center \(
    Y_{r}=\Max{f_{r}^{d}} \) (see \ref{comequiv}).
     \end{description}

     \item[C] For some index \( r \), depending on the basic object
     \( (W_{0},(J_{0},b),E_{0}) \), the equivariant sequence
     constructed in \textbf{B} is a resolution (i.~e.  \(
     \Sing(J_r,b)=\emptyset \)).
\end{description}

\begin{parrafo}
We refer the reader to \cite{BodnarSchicho2000_1} for an
implementation of the algorithm treated in
\cite[Theorem~7.13]{EncinasVillamayor2000}, to see how it works on
examples.

Condition {\bf C} says that for $i=0,1,\ldots,k$, the functions
\begin{equation*}
f_i:\Sing(J_i,b) \to I_{d}
\end{equation*}
define a resolution  of the basic object \((W_0,(J_0,b),E_0)\),
with
         centers \( \Max{f_{i}^{d}} \). We will refer to it as  {\em the
resolution defined by the algorithm}, or {\em the resolution
defined by the functions $f_i$.} Note that {\bf B} says that this
resolution is equivariant.

Let \( \mathcal{B}=(W_0,(J_0,b),E_0) \) be a basic object.  If \(
U_0\subset W_0\) is a non-empty open set, then we set the
\emph{restriction} of the basic object to be \(
(W_0,(J_0,b),E_0)_{U_0}=(U_0,(J|_{U_0},b),E_{U_0}) \), where \(
J|_{U_0} \) is the restriction of the sheaf of ideal to \( U_0\)
and \( E_{U_0}=\{H\cap U_0\mid H\in E\} \).

If
\begin{equation} \label{eqMainSeqResol1}
                  (W_{0},(J_{0},b),E_{0})\longleftarrow\cdots \longleftarrow
         (W_{N-1},(J_{N-1},b),E_{N-1})\longleftarrow
         (W_{N},(J_{N},b),E_{N})
              \end{equation}
is the resolution defined by the algorithm, it induces naturally a
sequence

\begin{equation} \label{eqMainSeqResol2}
                  (W_{0},(J_{0},b),E_{0})_{U_0}\longleftarrow\cdots
\longleftarrow
         (W_{N-1},(J_{N-1},b),E_{N-1})_{U_{N-1}}\longleftarrow
         (W_{N},(J_{N},b),E_{N})_{U_N}
              \end{equation}
where each $U_k$ is an open subset in $W_k$ (the pull back of
$U_0$ in $ W_k$).

Each function $f_i:\Sing(J_i,b) \to I_{d}$ induces by restriction
a function on $U_i \cap \Sing(J_i,b) $. Note that if \(
\Max{f_{i}^{d}}\cap U_i=\emptyset  \), then
$$(W_{i},(J_{i},b),E_{i})_{U_i}\longleftarrow
(W_{i+1},(J_{i+1},b),E_{i+1})_{U_{i+1}}$$ is the identity map, and
hence can be neglected from (\ref{eqMainSeqResol2}).
\end{parrafo}

\begin{parrafo}

\textbf{Properties of the algorithm.} Algorithmic principalization
has the following properties (see \cite[p. 192]{EncinasVillamayor2000}):
\begin{description}
     \item[p1] The functions defining the algorithmic resolution of the
     d-dimensional basic object $(W_{0},(J_{0},b),E_{0})_{U_0}$ are the
     restriction of the functions $f_k$ defining \ref{eqMainSeqResol1}.
     In particular \ref{eqMainSeqResol2} is the resolution of
     $(W_{0},(J_{0},b),E_{0})_{U_0}$ defined by the algorithm.

     \item[p2] For the resolution defined by the algorithm, say
     (\ref{eqMainSeqResol1}): \[ \max{f_{0}^{d}}>\max{f_{1}^{d}}>\cdots
     >\max{f_{N-1}^{d}} \]

     \item[p3] If \( J_{0} \) is the ideal of a regular pure
     dimensional subvariety \( X_{0} \), \( E_{0}=\emptyset \) and \(
     b=1 \), then the function \( f_{0}^{d} \) is constant.

     \item[p4] For any \( i=0,\ldots,N-1 \), the closed set \(
     \Max{f_{i}^{d}} \) is smooth and equidimensional.  Furthermore,
     the dimension is determined by the value \( \max{f_{i}^{d}} \).
     \end{description}
\end{parrafo}

\begin{parrafo}
It follows from Property {\bf p1)}, that if \( \xi\in
\Sing(J_{i},b) \), \( i=0,\ldots,r-1 \), and \( \xi\not\in Y_{i}
\), then \( f_{i}^{d}(\xi)=f_{i+1}^{d}(\xi') \) via the natural
identification of the point \( \xi \) with a point \( \xi' \) of
\( \Sing(J_{i+1},b) \)
\end{parrafo}

\begin{parrafo} \label{ProofMain}
\textbf{Proof of theorem~\ref{MainTh}}.
Fix notation as in theorem~\ref{MainTh}, and consider the basic
object
\begin{equation*}
(W_0,(J_0,1),E_0),
\end{equation*}
where $W_0=W$, $J_0=\mathcal{I}(X)$ and $E_0= \emptyset$. Clearly
$X= \Sing(J_0,1)$.

Take $U=W\setminus\Sing(X)$.  By {\bf p3)} we know that the function
$$f_0: \Sing(J_0,1)\to (I_d,\leq)$$
is constant on the restriction to
$ \Sing(J_0,1) \cap U$ ( on the restriction $(W_0,(J_0,1),E_0)_U$).
Let $a(d)$ denote this constant value along the points in $U \cap X$.

\medskip

By \ref{AlgResol} {\bf C)}, we know that the algorithm provides a
resolution of the basic object $(W_0,(J_0,1),E_0)$ by means of a
finite sequence of blow-ups
\begin{equation}
\label{pruebaclasica}
(W_0,(J_0,1),E_0)\longleftarrow(W_1,(J_1,1),E_1)
\longleftarrow\ldots\longleftarrow
(W_N,(J_N,1),E_N),
\end{equation}
at permissible centers $Y_i\subset \Sing(J_i,b)$ for
$i=0,1,\ldots,N-1$. Therefore $\Sing(J_N,b)= \emptyset $, and by
{\bf p1)} and {\bf p2)}, there must be an index
$k\in\{0,1,\ldots,N\}$ such that $\max f_k= a(d)$. Such index $k$
is unique by {\bf p2)}.

\medskip

Now  $U$ can be identified with an open set, say $U$ again, of
$W_k$ (note that the centers of the transformations in sequence
(\ref{pruebaclasica}) are defined by $\Max f_i$ and $\max f_i >
a(d)$ for $i<k$).  If $X_k$ denotes the strict transform of $X$ in
$W_k$, then
\begin{equation*}
        X_k \cap U=X \cap U= \Max f_k \cap U.
\end{equation*}
Since $X \cap U=\mbox{Reg}(X)$ is dense in $X$, it follows that
$X_k$ is the union of some of the components of $\Max f_k$, and
hence it is regular and has normal crossing with the exceptional
components by  Definition \ref{AlgResol} ({\bf A)}).  This proves
(i) and (ii) of Theorem~\ref{MainTh}.

\medskip

Now it only remains to show that the resolution of singularities
of $X$ that we have achieved is equivariant. This follows now from
\ref{MainEx}, together with the equivariant resolution of the
basic object $(W_0,(J_0,1),E_0)$ provided by the algorithm of
resolution of basic objects.

\end{parrafo}

\end{document}